\documentclass[12pt]{article}
\pdfoutput=1
\usepackage{amsfonts,hyperref}
\usepackage{color}
\usepackage{appendix}
\usepackage{amsmath}
\usepackage{graphicx}
\usepackage[latin1]{inputenc}
\usepackage[french,english]{babel}
\usepackage{geometry}
\usepackage{etoolbox}
\usepackage{fancyhdr}
\usepackage{url}
\usepackage{dsfont}
\usepackage{stmaryrd}
\patchcmd{\thebibliography}{\section*}{\section}{}{}

\addto\captionsfrench{}
\addto\captionsenglish{}

\newcommand\encadremath[1]{\vbox{\hrule\hbox{\vrule\kern8pt
\vbox{\kern8pt \hbox{$\displaystyle #1$}\kern8pt}
\kern8pt\vrule}\hrule}}
\def\enca#1{\vbox{\hrule\hbox{
\vrule\kern8pt\vbox{\kern8pt \hbox{$\displaystyle #1$}
\kern8pt} \kern8pt\vrule}\hrule}}

\newcommand\figureframex[3]{
\begin{figure}[bth]
\hrule\hbox{\vrule\kern8pt
\vbox{\kern8pt \vbox{
\begin{center}
{\mbox{\epsfxsize=#1.truecm\epsfbox{#2}}}
\end{center}
\caption{#3}
}\kern8pt}
\kern8pt\vrule}\hrule
\end{figure}
}

\makeatletter
\@addtoreset{equation}{section}
\makeatother
\newtheorem{theorem}{Theorem}[section]

\newtheorem{remark}{Remark}[section]
\newtheorem{proposition}{Proposition}[section]
\newtheorem{lemma}{Lemma}[section]
\newtheorem{corollary}{Corollary}[section]
\newtheorem{definition}{Definition}[section]

\def\br{\begin{remark}\rm\small}
\def\er{\end{remark}}
\def\bt{\begin{theorem}}
\def\et{\end{theorem}}
\def\bd{\begin{definition}}
\def\ed{\end{definition}}
\def\bp{\begin{proposition}}
\def\ep{\end{proposition}}
\def\bl{\begin{lemma}}
\def\el{\end{lemma}}
\def\bc{\begin{corollary}}
\def\ec{\end{corollary}}
\def\beaq{\begin{eqnarray}}
\def\eeaq{\end{eqnarray}}
\newcommand{\proof}[1]{{\noindent \bf proof:}\par
{#1} $\square$}

\newcommand{\beq}{\begin{equation}}
\newcommand{\eeq}{\end{equation}}
\newcommand{\bea}{\begin{eqnarray}}
\newcommand{\eea}{\end{eqnarray}}
\newcommand{\beqq}{\begin{equation*}}
\newcommand{\eeqq}{\end{equation*}}
\newcommand{\beaa}{\begin{eqnarray*}}
\newcommand{\eeaa}{\end{eqnarray*}}

\newcommand{\Pint}{{\int\kern -1.em -\kern-.25em}}

\renewcommand \P {{\mathbb{P}}}

\begin{document}

\sloppy

\pagestyle{empty}
\vspace{10pt}
\begin{center}
{\large \bf {Asymptotic expansion of Toeplitz determinants of an indicator function with discrete rotational symmetry and powers of random unitary matrices}}
\end{center}
\vspace{3pt}
\begin{center}
\textbf{O. Marchal$^\dagger$}
\end{center}
\vspace{15pt}

$^\dagger$ \textit{Universit\'{e} de Lyon, CNRS UMR 5208, Universit\'{e} Jean Monnet, Institut Camille Jordan, Institut Universitaire de France}
\footnote{olivier.marchal@univ-st-etienne.fr}

\vspace{30pt}

{\bf Abstract}: In this short article we propose a full large $N$ asymptotic expansion of the probability that the $m^{\text{th}}$ power of a random unitary matrix of size $N$ has all its eigenvalues in a given arc-interval centered in $1$ when $N$ is large. This corresponds to the asymptotic expansion of a Toeplitz determinant whose symbol is the indicator function of several intervals having a discrete rotational symmetry. This solves and improves a conjecture left opened by the author in \cite{ToeplitzMarchal2020}. It also provides a rare example of the explicit computation of a full asymptotic expansion of a genus $g>0$ classical spectral curve, including the oscillating non-perturbative terms, using the topological recursion.

\vspace{15pt}
\pagestyle{plain}
\setcounter{page}{1}


\section{Introduction and summary of the results}

Toeplitz determinants of size $N$ and their large $N$ asymptotic expansions are classical problems in probability in relation with the spectrum of random unitary matrices. Indeed, it is well-known that the Toeplitz determinant with symbol $f$ :
\beq \label{ToepDet} D_N(f)= \det \left(T_{p,q}=t_{p-q}(f)\right)_{1\leq p,q\leq N} \eeq
with discrete Fourier coefficients given by:
\beq \label{FourierCoeff} t_k(f)=\frac{1}{2\pi}\int_0^{2\pi}f(e^{i\theta})e^{-ik\theta} d\theta \,\,\,\,,\,\,\,\,  \forall\, k\in \mathbb{Z}\eeq
can be alternatively written as :
\beq \label{IntExp} D_N(f)=\frac{1}{(2\pi)^N\,N!}\int_{[-\pi,\pi]^N} d\theta_1\dots d\theta_N \left(\prod_{k=1}^N f(e^{i\theta_k})\right)  \prod_{1\leq p<q\leq N}\left|e^{i \theta_p}-e^{i\theta_q}\right|^2
\eeq
In particular, the $N$-dimensional integrand corresponds to the distribution of the eigenvalues of a random unitary matrix of size $N$ (i.e. the Circular Unitary Ensemble (CUE)) with an external potential characterized by the symbol $f$. Toeplitz determinants can be studied using orthogonal polynomials and Riemann-Hilbert problems (See \cite{FahsThesis} for a recent and accessible overview). They are also determinantal point processes and thus have strong relations with integrable systems (See \cite{MehtaBook2004} for a review) and Fredholm determinants (See \cite{BorodinOkounkovFredholmDetToeplitzDet2000} and references therein). Situations similar to the case studied in this article are  \cite{DeiftItsZhouRHToeplitz1997} where the authors considered the probability for the sine point process to have any number of gaps, while large gaps asymptotics in dimension two (involving random normal matrices) have recently been studied in \cite{GapsAnnuliRandomUnitaryMatricesCharlier2021}.

\medskip

Many results regarding Toeplitz determinants are available in the literature since the notion of Toeplitz determinants dates back to O. Toeplitz \cite{OttoToeplitz1907,OTTOToeplitz1911} in $1907$. The first important result in the case of regular and strictly positive symbols are those of Szeg\"{o} \cite{Szego1915} in $1915$. These results have been refined by many authors \cite{Szego1952,Ibragimov1968,GolinskiiIbragimov1971,Widom1971} up to Johansson work \cite{JohanssonStrongSzego1988} giving rise to the strong Szeg\"{o} theorem:
\beq \ln D_N(f)=\frac{N}{2\pi}\int_{0}^{2\pi} \ln f(e^{i\theta}) \, d\theta +\sum_{j=1}^\infty j(\ln f)_j (\ln f)_{-j} +o(1) \eeq
whenever $f$  satisfies $\underset{j=1}{\overset{\infty}{\sum}} |j| |(\ln f)_j|^2<+\infty$.

\medskip

However when the symbol is not regular enough or not strictly positive, the asymptotic expansion differs from the one above and much progress has been made during the past decades in the study of these cases. For example, Fisher-Hartwig singularities (i.e. isolated zeros or singular points on the unit circle) have been dealt with by many authors \cite{FisherHartwigSing1968,FisherHartwigBasorMorrison1994,DeifItsKrasovsky2011,DeifItsKrasovskyToepSingular2014,CharlierAsymptHankelMultiCutFHSingUpToConstantTerm2021,BlackStoneCharlierLenellsGapProbaBulkAiry2021} using orthogonal polynomials and Riemann-Hilbert problems.

\medskip

The case of a symbol equal to the indicator function of an interval, i.e. $f=\mathds{1}_{[\alpha,\beta]}$, corresponding to the gap probability in random unitary matrices, has been studied by H. Widom in \cite{Widom1974} and refined in \cite{WidomConstantGapProbaThirdProof2007,DuitsKozhanRelativeSzegoAsymptToep2017}. Finally, the full asymptotic expansion was presented in \cite{ToeplitzMarchal2020} (Theorem $2.3$). Coefficients of the asymptotic expansion are directly related to the output of the topological recursion defined by B. Eynard and N. Orantin \cite{EynardOrantin2007} applied to the genus $0$ classical spectral curve $y^2=\frac{1}{\cos^2\left(\frac{|\beta-\alpha|}{4}\right)(1+x^2)^2\left(x^2-\tan^2\left(\frac{|\beta-\alpha|}{4}\right)\right)}$. Symbols supported on a single interval with singular points inside the support have recently been studied in \cite{KrasovskyToeplitzSingleArcWithZerosInside2006}.

\medskip

Finally, in the past years, progress has been made in the case of a symbol supported on a union of intervals (known as the multi-cut case). The theoretical form of the asymptotic expansion was proved in \cite{BorotGuionnetBetaMultiCut2013,BorotGuionnetKozBetaMultiCut2015} generalizing similar results available for the one-cut case \cite{BorotGuionnetBetaOneCut2013}. The first leading coefficients have been obtained in the case of two intervals \cite{CharlierClayesToepSymbolWithAGap2015,BesselDeterminant2021}. Leading terms, including the constant terms, have been obtained in Airy and Sine kernels determinants on two large intervals in respectively \cite{KrasovskyMaroudas2021} and \cite{FahsKrasovsky2020}. Oscillatory terms in the case of the Airy kernel have also been computed in \cite{BlackStoneCharlierLenellsOsciAiryKernelTwoIntervals2020}. When the symbol is an indicator function on a union of intervals, i.e. $f=\mathds{1}_{\underset{k=1}{\overset{m}{\bigcup}}[\alpha_k,\beta_k]}$, only the first two leading terms are known. In particular, a general formula for the constant term is still missing and very few is known. In \cite{ToeplitzMarchal2020}, a conjecture in the case when the symbol has an additional discrete rotational symmetry was proposed for the constant term.
\bigskip 

\bigskip

The main result of this article is to prove the conjecture left opened in \cite{ToeplitzMarchal2020} and propose a full large $N$ asymptotic expansion of the Toeplitz determinant $D_N(m,\epsilon) $ with symbol $f_{m,\epsilon}=\mathds{1}_{\underset{k=0}{\overset{m-1}{\bigcup}}\left[\frac{2\pi k}{m}-\frac{\pi \epsilon}{m}, \frac{2\pi k}{m}+\frac{\pi \epsilon}{m}\right]}$ in Theorem \ref{MainTheo}. The proof is based on a combination of a result of B. Fahs \cite{FahsThesis} factorizing the initial Toeplitz determinant and results of \cite{ToeplitzMarchal2020} in the one interval case. As a byproduct, we obtain the full large $N$ asymptotic expansion of the probability that the $m^{\text{th}}$ power of a random unitary matrix of size $N$ denoted $U_N$ (i.e. a random Haar distributed $N\times N$ unitary matrix) has all its eigenvalues closed to $1$ (See Corollary \ref{CorollaryUnitaryMatrices}). In the end, we may sum up our main results by:
\bea\label{Dnf} D_{N}(f_{m,\epsilon})&=&\P\left( ||U_N^{\,m}-I_N||_2 \leq \sqrt{2}\sin\frac{\pi \epsilon}{2}\right)\cr
&=&\P\left( \text{All eigenvalues of }\, U_N^{\,m} \, \text{ belong to }\, \left\{e^{it}\,,\, t\in [-\pi \epsilon,\pi \epsilon]\right\}\right) \cr
&\overset{N=n_1m+n_2}{\underset{n_1\to \infty}{=}}&n_1^2m \ln \left(\sin\frac{\pi \epsilon}{2}\right) +2n_1 n_2 \ln \left(\sin\frac{\pi \epsilon}{2}\right)-\frac{m}{4}\ln n_1\cr
&&+m\left(\frac{1}{12}\ln 2-3\,\zeta'(-1)\right)-\frac{m}{4}\ln \left(\cos\frac{\pi \epsilon}{2}\right)+ n_2\ln \left(\sin\frac{\pi \epsilon}{2}\right) -\frac{n_2}{4n_1}\cr
&&-\sum_{l=1}^{\infty} \left(mF_{\text{Top. Rec.}}^{(l+1)}(\epsilon) -\frac{n_2}{8l}+n_2\sum_{j=1}^{l-1} \binom{2l-1}{2j-1} F_{\text{Top. Rec.}}^{(j+1)}(\epsilon)\right)\frac{1}{n_1^{2l}}\cr
&&+n_2\sum_{l=1}^{\infty}\left(-\frac{1}{4(2l+1)} +\sum_{j=1}^l \binom{2l}{2j-1} F_{\text{Top. Rec.}}^{(j+1)}(\epsilon)\right)\frac{1}{n_1^{2l+1}}+ o(n_1^{-\infty})\cr
&&
\eea
where $||.||_2$ is the Euclidean norm, $N=n_1m+n_2$ is the Euclidean division of $N$ by $m$ and $\left(F_{\text{Top. Rec.}}^{(k)}(\epsilon)\right)_{k\geq 0}$ are the free energies associated to the genus $0$ classical spectral curve 
\beqq y^2=\frac{1}{\cos^2\left(\frac{\pi \epsilon}{2}\right)(1+x^2)^2\left(x^2-\tan^2\left(\frac{\pi \epsilon}{2}\right)\right)}\eeqq
computed by the topological recursion of \cite{EynardOrantin2007} (See \eqref{FgTopRec} for the first free energies that are also presented in full details in \cite{ToeplitzMarchal2020}). Note also that the first terms of \eqref{Dnf}, up to $O\left(N^{-1}\right)$, where obtained by B. Fahs in \cite{FahsThesis}. 

\medskip

This complements existing results on powers of random unitary matrices \cite{DiaconisShahshahani1994,DuitsJohanssonPowerUnitaryToeplitz2010,KrasovskyAspectsToeplitzDet2011,UnitaryMarchal2014,CharlierClaeysThinningCUE2017} and might lead to sharper control when no symmetry arise. The result is also interesting for the topological recursion community. Indeed, the classical spectral curve associated to the symbol $f_{m,\epsilon}$ is of genus $g=m-1$ (See Proposition \ref{TopRecGenus}). Obtaining an explicit expression for the full asymptotic expansion of the partition function of a classical spectral curve of strictly positive genus is known to be a very difficult problem for which no other example are currently known. Thus, the present article provides the first example of an arbitrary genus spectral curve for which the full asymptotic expansion, including the non-perturbative oscillatory terms can be computed up to any order by the means of topological recursion. In particular, it opens the way to similar cases where discrete symmetries shall be used efficiently to reduce an impractical strictly positive spectral curve problem to a mere implementation of the topological recursion on a genus zero classical spectral curve.   
 
\section{Asymptotic expansion of Toeplitz determinants of an indicator function with discrete rotational symmetry} 

\subsection{Definitions and existing results}

\begin{definition}\label{DefinitionToeplitz} Let $m\geq 1$, $\epsilon\in (0,1)$ and $N\geq 1$. Let $\mathcal{C}$ be the unit circle. We define:
\beqq I_m(\epsilon)=\underset{k=0}{\overset{m-1}{\bigcup}}\left[\frac{2\pi k}{m}-\frac{\pi \epsilon}{m}, \frac{2\pi k}{m}+\frac{\pi \epsilon}{m}\right]\,\,\,,\,\, \mathcal{C}_m(\epsilon)=\{e^{i\theta}, \theta\in I_m(\epsilon)\}\subset \mathcal{C}\eeqq
In other words, $\mathcal{C}_m(\epsilon)$ corresponds to the union of $m$ arc-intervals of size $\frac{2\pi \epsilon}{m}$ centered at $\left(\gamma_k=e^{\frac{2\pi  i k}{m}}\right)_{k=0}^{m-1}$. The corresponding Toeplitz determinant is defined by:
\beaa D_N(m,\epsilon)&=&\frac{1}{(2\pi)^N\,N!}\int_{I_m(\epsilon)^{\, N}} d\theta_1\dots d\theta_N  \prod_{1\leq p<q\leq N}\left|e^{i \theta_p}-e^{i\theta_q}\right|^2\cr
&=&D_N(f_{m,\epsilon}) \text{ with } f_{m,\epsilon}=\mathds{1}_{I_m(\epsilon)}\cr
&=& \det \left(T_{p,q}=t_{p-q}(m,\epsilon)\right)_{1\leq p,q\leq N}
\eeaa
with  $t_0(m,\epsilon)=\epsilon$ and  $t_k(m,\epsilon)= \epsilon \sin_c\left(\frac{\pi \epsilon k}{m}\right) \delta_{k\,\equiv\, 0\,[m]} \,\text{ for } k\neq 0$. Note in particular that $D_N(1,m\epsilon)=D_N(\mathds{1}_{[-\pi \epsilon,\pi \epsilon]})$ is independent of $m$. By convention, we also set $D_0(m,\epsilon)=1$ for any values of $m$ and $\epsilon$.\footnote{Note that the definition of the Toeplitz determinant $D_N(m,\epsilon)$ remains unchanged by any global rotation of the angles (i.e. $\theta_i\mapsto \theta_i+\alpha$ for all $i\in \llbracket 1,N\rrbracket$). In particular, in \cite{ToeplitzMarchal2020}, for $m$ even, all angles where shifted by $\frac{2\pi}{2m}$ so that the function $\theta\mapsto \tan \frac{\theta}{2}$ could be applied on all intervals defining $I_m(\epsilon)$.} 
\end{definition}

Note that $\mathcal{C}_m(\epsilon)$ is invariant under the rotation $z\mapsto ze^{\frac{2\pi i}{m}}$. It may be represented graphically: 

\begin{center}
\includegraphics[width=12cm]{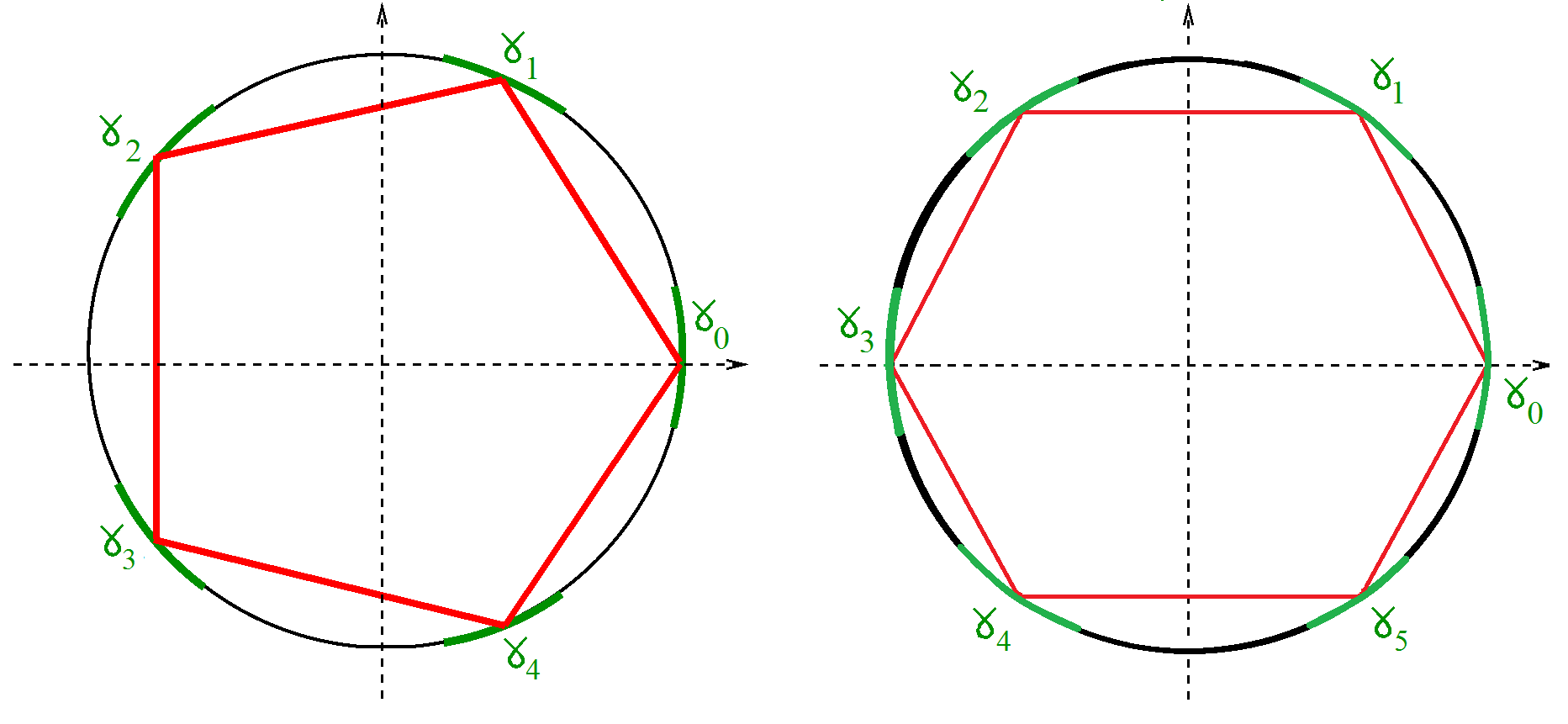}

\textit{Fig. $1$: Illustration (in green) of the set $\mathcal{C}_m(\epsilon)$ for $m=5$ (left) or $m=6$ (right).}
\end{center}

\medskip

The correspondence between the Toeplitz determinant and the $N$-dimensional integral is straightforward and details can be found in \cite{ToeplitzMarchal2020} with additional reformulations of the problem. Note also that the determinant formulation (with a Toeplitz matrix having lots of vanishing entries) is particularly convenient for numerical simulations. In \cite{ToeplitzMarchal2020}, it was proved using results of \cite{BorotGuionnetKozBetaMultiCut2015} that the Toeplitz determinant $D_N(m,\epsilon)$ has a full large $N$ asymptotic expansion related to the topological recursion of \cite{EynardOrantin2007}:

\begin{proposition}\label{TopRecGenus}The Toeplitz determinant $D_N(m,\epsilon)$ has the following large $N$ asymptotic expansion:
\small{\bea \label{GeneralForm2R1}&&D_N(m,\epsilon)\overset{N\to \infty}{=}\frac{N^{N+\frac{1}{4}m}}{N!} \exp \left(\sum_{k=-2}^\infty N^{-k}F_{\boldsymbol{\epsilon}^\star}^{\{k\}}\right)\cr
&&\left\{\sum_{j\geq 0}\,\sum_{\begin{subarray}{l}l_1,\dots,l_j\geq 1\\ k_1,\dots,k_j\geq -2 \\ \underset{i=1}{\overset{j}{\sum}} (l_i+k_i)>0\end{subarray}} \frac{N^{-\left(\underset{i=1}{\overset{j}{\sum}} (l_i+k_i)\right)}}{j!}
\left(\underset{i=1}{\overset{j}{\bigotimes}} \frac{ F^{\{k_i\},(l_i)}_{\boldsymbol{\epsilon}^\star}}{l_i!}\right) \cdot \nabla_\nu^{\otimes\left(\underset{i=1}{\overset{j}{\sum}} l_i\right)} \right\}\Theta_{-N\boldsymbol{\epsilon}^\star}\left(\mathbf{0}\big|F^{\{-2\},(2)}_{\boldsymbol{\epsilon}^\star}\right) + o(N^{-\infty})\cr
&&\eea}
\normalsize{where} $\boldsymbol{\epsilon}^\star=\left(\frac{1}{m},\dots,\frac{1}{m}\right)^t$.
The coefficients $\left(F_{\boldsymbol{\epsilon}^\star}^{\{k\}}\right)_{k\geq -2}$ are related to the free energies $\left(F^{(g)}\right)_{g\geq 0}$ computed from the topological recursion applied to the classical spectral curve: 
\begin{itemize} \item \normalsize{If} $m=2r+1$ is odd:
\small{\beq\label{SpecCurveGenus} y^2=\lambda_r\frac{\underset{j=-r}{\overset{r-1}{\prod}} \left(x-\tan\left(\frac{j\pi}{2r+1}+\frac{\pi}{2(2r+1)}\right)\right)^2}{(1+x^2)^2\underset{k=-r}{\overset{r}{\prod}}\left(x^2-\tan^2\left(\frac{\pi k}{2r+1}+\frac{\pi \epsilon}{2(2r+1)}\right)\right)} \text{ with } \lambda_r=\frac{\underset{k=0}{\overset{r-1}{\prod}} \cos^4\left(\frac{\pi k}{2r+1}+\frac{\pi}{2(2r+1)}\right)} {\underset{k=-r}{\overset{r}{\prod}}\cos^2\left(\frac{\pi k}{2r+1}+\frac{\pi \epsilon}{2(2r+1)}\right)}
\eeq}
\item \normalsize{If} $m=2s$ is even:
\small{\beq y^2= \frac{ x^2\underset{k=1}{\overset{s-1}{\prod}}\left(x^2-\tan^2\left(\frac{\pi k}{2s}\right)\right)^2}{(1+x^2)^2\underset{k=-(s-1)}{\overset{s}{\prod}}\left(x^2-\tan^2\left(\frac{\pi\left(k-\frac{1}{2}\right)}{2s}+\frac{\pi \epsilon}{4s}\right)\right)} \text{ with } \lambda_s= \frac{\underset{k=1}{\overset{s-1}{\prod}} \cos^2\left(\frac{\pi\left(k-\frac{1}{2}\right)}{2s}\right)}{\underset{k=-(s-1)}{\overset{s}{\prod}} \cos^2\left(\frac{\pi\left(k-\frac{1}{2}\right)}{2s}-\frac{\pi\epsilon}{4s}\right)}
\eeq}\normalsize{}
\end{itemize}
by:
\beaa \forall \,k\geq -1\,&:&\,\, F_{\boldsymbol{\epsilon}^\star}^{\{2k\}}=-F^{(2k+2)}+f_{2k} \text{ with } f_{2k} \text{ independent of } \epsilon \cr
 \forall \,k\geq -1\,&:&\,\, F_{\boldsymbol{\epsilon}^\star}^{\{2k+1\}}=f_{2k+1} \text{ with } f_{2k+1} \text{ independent of } \epsilon
 \eeaa
The expression of $\left(f_{j}\right)_{j\geq 1}$ can be obtained from the limit $\epsilon\to 0$ of $D_N(m,\epsilon)$ given by equation $3.34$ of \cite{ToeplitzMarchal2020}.
\end{proposition}

In the former proposition, $\Theta$ is the Siegel theta function:
\beqq \Theta_{\boldsymbol{\gamma}}(\boldsymbol{\nu},\mathbf{T})=\sum_{\mathbf{k}\in \mathbb{Z}^g} \text{exp}\left(-\frac{1}{2}(\mathbf{k}+\boldsymbol{\gamma})^t\,\cdot \mathbf{T}\cdot(\mathbf{k}+\boldsymbol{\gamma})+\boldsymbol{\nu}^t\,\cdot (\mathbf{k}+\boldsymbol{\gamma})\right)\eeqq
and $F^{\{2k\},(l)}_{\boldsymbol{\epsilon}}$ are defined as the $l^{\text{th}}$ derivative of the coefficient $F_{\boldsymbol{\epsilon}}^{\{2k\}}$ relatively to the filling fractions $\boldsymbol{\epsilon}=(\epsilon_1,\dots,\epsilon_{m})^t \in\{\mathbf{u}\in (\mathbb{Q}_+)^{m}\,/ \, \underset{i=1}{\overset{m}{\sum}}u_i=1\}$. We stress here that large $N$ asymptotic expansions presented in this article (denoted with $o(N^{-\infty})$ as in \cite{BorotGuionnetKozBetaMultiCut2015})  are to be understood as asymptotic expansions up to any arbitrary large negative power of $N$ as defined in \cite{BorotGuionnetKozBetaMultiCut2015}.

\begin{remark} In the former proposition, quantities $\left(F^{(g)}\right)_{g\geq 0}$ are computed by the standard topological recursion of \cite{EynardOrantin2007} applied to the genus $0$ classical spectral curve \eqref{SpecCurveGenus}. For example, the leading orders $F^{(0)}$ and $F^{(1)}$ are given by
\beqq F^{(0)}=\frac{1}{m}\ln\left(\sin \frac{\pi \epsilon}{2}\right)\,\,,\,\, F^{(1)}=-\frac{m}{4}\ln\left(\cos \frac{\pi \epsilon}{2}\right)
\eeqq
\end{remark}

\begin{remark}The expression 
\beqq \left\{\sum_{j\geq 0}\,\sum_{\begin{subarray}{l}l_1,\dots,l_j\geq 1\\ k_1,\dots,k_j\geq -2 \\ \underset{i=1}{\overset{j}{\sum}} (l_i+k_i)>0\end{subarray}} \frac{N^{-\left(\underset{i=1}{\overset{j}{\sum}} (l_i+k_i)\right)}}{j!}
\left(\underset{i=1}{\overset{j}{\bigotimes}} \frac{ F^{\{k_i\},(l_i)}_{\boldsymbol{\epsilon}^\star}}{l_i!}\right) \cdot \nabla_\nu^{\otimes\left(\underset{i=1}{\overset{j}{\sum}} l_i\right)} \right\}\Theta_{-N\boldsymbol{\epsilon}^\star}\left(\mathbf{0}\big|F^{\{-2\},(2)}_{\boldsymbol{\epsilon}^\star}\right)
\eeqq
in \eqref{GeneralForm2R1} is presented in details in \cite{BorotGuionnetBetaMultiCut2013}. In particular equations $1.29$, $1.31$ and $1.32$ of \cite{BorotGuionnetBetaMultiCut2013} provide the first two orders and details. Since the expressions are rather long and of no particular interest for the present work, we choose not to reproduce them here.
\end{remark}

As one can see, the asymptotic expansion of $D_N(m,\epsilon)$ presented in Proposition \ref{TopRecGenus} is rather complicated. In \cite{ToeplitzMarchal2020} the discrete rotational symmetry of the problem was used to obtain an exact expression for the classical spectral curve and to obtain the optimal (uniform) filling fraction vector $\boldsymbol{\epsilon}^\star$. Unfortunately, no simplification of the Siegel Theta function could be made. In practice, only the first two leading coefficients are easily computable using Proposition \ref{TopRecGenus} and one only gets:
\beq\label{FirstTwoOrders} \ln D_N(m,\epsilon)= \frac{N^2}{m}\ln \left(\sin\left(\frac{\pi \epsilon}{2}\right)\right)-\frac{m}{4}\ln N+O(1)\eeq

\subsection{A full large $N$ asymptotic expansion of $D_N(m,\epsilon)$} 

In \cite{FahsThesis}, B. Fahs observed that the rotational invariance can be used in a much stronger way. Indeed, at the level of Toeplitz determinants we have:

\begin{proposition}[Factorization of $D_N(m,\epsilon)$ (Proposition 1.1.3 of \cite{FahsThesis})]\label{Factorization} Let $N\geq 1$ and write $N=n_1 m+n_2$ with $n_2\in\llbracket 0, m-1\rrbracket$ the Euclidean division of $N$ by $m$\footnote{Note that $n_1$ and $n_2$ depends on $N$ ($n_2$ a periodic function of $N$ with period $m$ while $n_1=\left\lfloor \frac{N}{m}\right\rfloor$) but we shall not write down the dependence explicitly for clarity.}. Then
\beqq D_N(m,\epsilon)= D_{n_1}(1,m\epsilon)^{m-n_2} D_{n_1+1}(1,m\epsilon)^{n_2}=D_{n_1}(\mathds{1}_{[-\pi \epsilon,\pi \epsilon]})^{m-n_2} D_{n_1+1}(\mathds{1}_{[-\pi \epsilon,\pi \epsilon]})^{n_2}\eeqq
\end{proposition}

In other words the Toeplitz determinant $D_N(m,\epsilon)$ on $m$ intervals can be expressed only in terms of the Toeplitz determinant with symbol given by the indicator function on a single interval $I=[-\pi \epsilon,\pi \epsilon]$.  It turns out that the large $N$ asymptotic expansion of the one interval case is completely known:

\begin{proposition}[Asymptotic expansion of $D_n(\mathds{1}_{[-\pi \epsilon,\pi \epsilon]})$ (Theorem $2.3$ of \cite{ToeplitzMarchal2020})]\label{OneCutCase} For $m\geq 1$ and $\epsilon\in(0,1)$ we have:
\bea \ln D_n(\mathds{1}_{[-\pi \epsilon,\pi \epsilon]})&\overset{n\to \infty}{=}&n^2\ln\left(\sin\left(\frac{\pi \epsilon}{2}\right)\right)-\frac{1}{4}\ln n -\frac{1}{4}\ln\left(\cos\left(\frac{\pi \epsilon}{2}\right)\right)\cr 
&&+3\,\zeta'(-1)+\frac{1}{12}\ln 2- \sum_{g=2}^\infty F_{\text{Top. Rec.}}^{(g)}(\epsilon) n^{2-2g}+ o(n^{-\infty})\eea
where the coefficients $\left(F_{\text{Top. Rec.}}^{(g)}(\epsilon)\right)_{g\geq 2}$ are the free energies associated to the genus $0$ classical spectral curve 
\beq \label{SpecCurveGenus0}y^2=\frac{1}{\cos^2\left(\frac{\pi \epsilon}{2}\right)(1+x^2)^2\left(x^2-\tan^2\left(\frac{\pi \epsilon}{2}\right)\right)}\eeq
\end{proposition} 

Combining Propositions \eqref{Factorization} and \eqref{OneCutCase} allows for the computation of the full large $N$ asymptotic expansion of  $D_N(m,\epsilon)$.  Using the identity:
\beqq \frac{1}{(1+x)^l}=1+\sum_{k=1}^{\infty} \binom{k+l-1}{l-1} (-1)^k x^k \,\,\, ,\,\,\, \forall\, l\geq 1 \,\, \text{ and } |x|<1\eeqq
a straightforward but tedious computation gives our main theorem:

\begin{theorem}[Asymptotic expansion of $D_N(m,\epsilon)$]\label{MainTheo} Let $m\geq 1$ and write $N=n_1 m+n_2$ with $n_2\in\llbracket 0, m-1\rrbracket$ the Euclidean division of $N$ by $m$. Then
\bea\label{MainTheoAsymp} \ln D_{N}(m,\epsilon)&\overset{N=mn_1+n_2}{\underset{n_1\to \infty}{=}}&n_1^2m \ln \left(\sin\frac{\pi \epsilon}{2}\right) +2n_1 n_2 \ln \left(\sin\frac{\pi \epsilon}{2}\right)-\frac{m}{4}\ln n_1\cr
&&+m\left(\frac{1}{12}\ln 2-3\,\zeta'(-1)\right)-\frac{m}{4}\ln \left(\cos\frac{\pi \epsilon}{2}\right)+ n_2\ln \left(\sin\frac{\pi \epsilon}{2}\right) -\frac{n_2}{4n_1}\cr
&&-\sum_{l=1}^{\infty} \left(mF_{\text{Top. Rec.}}^{(l+1)}(\epsilon) -\frac{n_2}{8l}+n_2\sum_{j=1}^{l-1} \binom{2l-1}{2j-1} F_{\text{Top. Rec.}}^{(j+1)}(\epsilon)\right)\frac{1}{n_1^{2l}}\cr
&&+n_2\sum_{l=1}^{\infty}\left(-\frac{1}{4(2l+1)} +\sum_{j=1}^l \binom{2l}{2j-1} F_{\text{Top. Rec.}}^{(j+1)}(\epsilon)\right)\frac{1}{n_1^{2l+1}}+ o(n_1^{-\infty})\cr
&&
\eea
where $\left(F_{\text{Top. Rec.}}^{(k)}(\epsilon)\right)_{k\geq 0}$ are the free energies associated to the genus $0$ classical spectral curve 
\beqq y^2=\frac{1}{\cos^2\left(\frac{\pi \epsilon}{2}\right)(1+x^2)^2\left(x^2-\tan^2\left(\frac{\pi \epsilon}{2}\right)\right)}\eeqq
computed by the topological recursion of \cite{EynardOrantin2007}. For example, the first free energies are given by:
\bea \label{FgTopRec} F_{\text{Top. Rec.}}^{(0)}(\epsilon)&=&\ln 2-\ln\left(\sin\frac{\pi \epsilon}{2}\right)  \cr
F_{\text{Top. Rec.}}^{(1)}(\epsilon)&=&\frac{1}{4}\ln\left(\cos\left(\frac{\pi \epsilon}{2}\right)\right) \cr
F_{\text{Top. Rec.}}^{(2)}(\epsilon)&=&\frac{1}{64}-\frac{1}{32}\tan^2\left(\frac{\pi \epsilon}{2}\right)\cr
F_{\text{Top. Rec.}}^{(3)}(\epsilon)&=&-\frac{1}{256}-\frac{1}{128}\tan^2\left(\frac{\pi \epsilon}{2}\right)-\frac{5}{128}\tan^4\left(\frac{\pi \epsilon}{2}\right)
\eea
\end{theorem}

\begin{remark} Computations of $F_{\text{Top. Rec.}}^{(0)}(\epsilon)$ and $F_{\text{Top. Rec.}}^{(1)}(\epsilon)$ in \eqref{FgTopRec} are presented in details in Appendix $A$ of \cite{ToeplitzMarchal2020}. Computations of $F_{\text{Top. Rec.}}^{(2)}(\epsilon)$ and $F_{\text{Top. Rec.}}^{(3)}(\epsilon)$ require the computation of the first orders of the topological recursion of \cite{EynardOrantin2007}. Intermediate steps are provided in Appendix \ref{AppendixComputationTopRec} for completeness.
\end{remark}

Using the determinant formulation, we may test the large $N$ asymptotic expansion with the finite $N$ numerical computations or the Toeplitz determinant:

\begin{center}
\includegraphics[width=8cm]{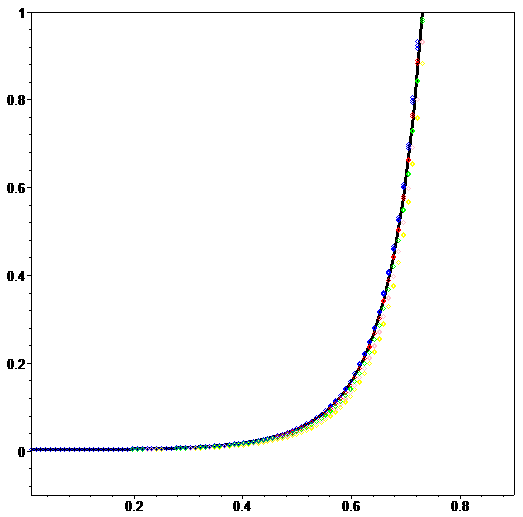}

\textit{Fig. $2$: Numerical computations of $\epsilon \mapsto \Big(D_N(m,\epsilon)- \big[n_1^2m \ln \left(\sin\frac{\pi \epsilon}{2}\right) +2n_1 n_2 \ln \left(\sin\frac{\pi \epsilon}{2}\right)-\frac{m}{4}\ln n_1
+m\left(\frac{1}{12}\ln 2-3\,\zeta'(-1)\right)-\frac{m}{4}\ln \left(\cos\frac{\pi \epsilon}{2}\right)+ n_2\ln \left(\sin\frac{\pi \epsilon}{2}\right) -\frac{n_2}{4n_1}+\frac{n_2}{8n_1^2}-\frac{m}{n_1^2}F_{\text{Top. Rec.}}^{(2)}(\epsilon)-\frac{n_2}{12 n_1^3}+\frac{2n_2}{n_1^3} F_{\text{Top. Rec.}}^{(2)}(\epsilon)+\frac{n_2}{16 n_1^4}-\frac{3n_2}{n_1^4}F_{\text{Top. Rec.}}^{(2)}(\epsilon) \big]\Big)\frac{n_1^4}{m}$ for $m=5$, $N\in\llbracket 90, 100\rrbracket$ (colored dots). The black curve is $\epsilon\mapsto F_{\text{Top. Rec.}}^{(3)}(\epsilon)=-\frac{1}{256}-\frac{1}{128}\tan^2\left(\frac{\pi \epsilon}{2}\right)-\frac{5}{128}\tan^4\left(\frac{\pi \epsilon}{2}\right)$. Colors of the points correspond to the different values of $n_2$ ($\{0, 1, 2, 3, 4\}$ are respectively matched with $\{$blue, red, green, pink, yellow$\}$) }
\end{center}

\begin{remark}Note that the asymptotic expansion \eqref{MainTheoAsymp} recovers \eqref{FirstTwoOrders}. Indeed, since $n_1=\left\lfloor \frac{N}{m}\right\rfloor$ and $n_2=N-m\left\lfloor \frac{N}{m}\right\rfloor$, the asymptotic expansion \eqref{MainTheoAsymp} can also be written as:
\bea\label{MainTheoAsymp2} &&\ln D_N(m,\epsilon)\overset{N\to \infty}{=}N^2 \ln \left(\sin\frac{\pi \epsilon}{2}\right) -\frac{m}{4}\ln N+ \frac{m}{4}\ln m+m\left(\frac{1}{12}\ln 2-3\,\zeta'(-1)\right)\cr
&&-\frac{m}{4}\ln \left(\cos\frac{\pi \epsilon}{2}\right)- \left(N-m\left\lfloor \frac{N}{m}\right\rfloor\right)\left(N-1-m\left\lfloor \frac{N}{m}\right\rfloor\right)\ln \left(\sin\frac{\pi \epsilon}{2}\right) 
+o(1)
\eea
In particular, this form provides the constant term of the large $N$ asymptotic expansion of $\ln D_N(m,\epsilon)$. This constant term was already obtained prior to this work by B. Fahs (Cf. eq. $1.61$ of \cite{FahsThesis}). It also verifies Conjecture $3.1$ of \cite{ToeplitzMarchal2020}. 
\end{remark}

\begin{remark}The asymptotic expansion presented in Theorem \ref{MainTheo} is performed in terms of  large $n_1$ and fixed $n_2$. However, one may easily obtain an asymptotic expansion in terms of $N$ and $n_2$ using the fact that $N=n_1 m+n_2$. Indeed, one can get the asymptotic expansion expressed in terms of  $N$ and $n_2$ by simply replacing $n_1$ into $\frac{N}{m}\left(1-\frac{n_2}{N}\right)$ in \eqref{MainTheoAsymp}, expand all corresponding quantities at large $N$ and reorder the series expansion. For example, powers like $n_1^{-l}$ provide $\underset{j=0}{\overset{\infty}{\sum}} \binom{j+l-1}{l-1}\frac{n_2^j m^l}{N^{j+l}}$ and $\ln n_1$ provides $\ln N-\ln m+\underset{k=1}{\overset{\infty}{\sum}} \frac{n_2^k}{k N^k}$. Since the final expression is more complicated than \eqref{MainTheoAsymp}, we do not see sufficient interest to write it down (the first orders being given by \eqref{MainTheoAsymp2}).  
\end{remark}

\section{Corollary for unitary matrices and discussion}
\subsection{Application to powers of random unitary matrices}

The connection between the Toeplitz determinant \eqref{IntExp} and the eigenvalues of a random unitary matrix of size $N$ provides the following corollary:

\begin{corollary}\label{CorollaryUnitaryMatrices}Let $U_N$ be a random unitary matrix of size $N$ drawn from the Haar measure on the set of unitary matrices and let $\epsilon\in (0,1)$. Let $m\geq 1$ be a given integer and write $N=n_1 m+n_2$ the Euclidean division of $N$ by $m$. Then the probability
\beqq p_{N,m}(\epsilon)= \P\left( ||U_N^{\,m}-I_N||_2 \leq \sqrt{2}\sin\frac{\pi \epsilon}{2}\right)=\P\left( \text{All eigenvalues of }\, U_N^{\,m} \, \text{ belong to }\, \mathcal{C}_1(\epsilon)\right) 
\eeqq
equals $D_N(m,\epsilon)$ and thus has a full large $N$ asymptotic expansion given by Theorem \ref{MainTheo}.
\end{corollary}

\proof{
Let $\left(u_1=e^{i\theta_1},\dots,u_N=e^{i\theta_N}\right)$ be the eigenvalues of $U_N^{\, m}$ with $(\theta_1,\dots,\theta_N)\in [0,2\pi]^N$. For a diagonalizable matrix $A$ of size $N$, it is well-known that $||A||_2$ corresponds to $\underset{j\in \llbracket 1,N\rrbracket}{\text{Max}} \{|\lambda_j|\,,\,  \lambda_j\, \text{ eigenvalues of } A\}$. Thus, $||U_N^{\, m}-I_N||_2=\underset{j\in \llbracket 1,N\rrbracket}{\text{Max}}\left\{\sqrt{2}\left|\sin \left(\frac{\theta_j}{2}\right)\right|\right\}$. For $\theta\in [0,2\pi]$, the inequality $\sqrt{2}\left|\sin\left( \frac{\theta}{2}\right)\right|\leq \sqrt{2}\sin\left(\frac{\pi \epsilon}{2}\right)$ is equivalent to $ 0\leq \sin \left(\frac{\theta}{2}\right)\leq \sin\left(\frac{\pi \epsilon}{2}\right)$ because $ \sin \frac{\theta}{2}$ is always non-negative. This inequality is equivalent to $\theta\in [0,\pi \epsilon]$ or $\frac{\theta}{2}\in \left[\pi-\frac{\pi \epsilon}{2},\pi\right]$, i.e. $\theta\in [0,\pi \epsilon]\cup[2\pi-\pi \epsilon, 2\pi]$. Hence, $||U_N^{\, m}-I_N||_2\leq \sqrt{2}\sin\left(\frac{\pi \epsilon}{2}\right)$ is equivalent to say that for all $j\in \llbracket 1,N\rrbracket$: $e^{i \theta_j}\in \mathcal{C}_1(\epsilon)$. Finally, this is equivalent to say that all eigenvalues $\left(e^{\frac{i\theta_j}{m}}\right)_{1\leq j \leq N }$ of $U_N$ belong to $\mathcal{C}_m(\epsilon)$ whose probability is given by $D_N(m,\epsilon)$.
}

\medskip

\medskip

We observe that Proposition \ref{Factorization} is compatible with the well-known results that for $m > N$, the eigenvalues of $U_N^{\, m}$ are i.i.d. random variables distributed uniformly on the unit circle \cite{DiaconisShahshahani1994}. Indeed, in Proposition \ref{Factorization}, this corresponds to $n_1=0$ and $n_2=N$, so that the factorization provides:
\beq \forall \,m> N\,\,:\,\, D_N(m,\epsilon)= D_{0}(1,m\epsilon)^{m-N} D_{1}(1,m\epsilon)^{N}=\left(\int_{-\pi\epsilon}^{\pi\epsilon} \frac{d\theta}{2\pi}\right)^N=\epsilon^N\eeq
Note that this case is not covered by the large $N$ expansion given by Theorem \ref{MainTheo}.

\subsection{Interest for the topological recursion formalism}

Theorem \ref{MainTheo} is interesting from the topological recursion perspective. Indeed, the classical spectral curve associated to the problem is given by \eqref{SpecCurveGenus} and the expected asymptotic expansion for such a spectral curve is given in Proposition \ref{TopRecGenus}. In general, the presence of a Theta function makes the expansion very difficult to compute in practice when the genus of the spectral curve is not $0$. In our case, it turns out that the discrete rotational symmetry of the problem allows for the explicit computation because of the factorization \ref{Factorization}. This implies that the Eynard-Orantin differentials and free energies associated to the genus $m-1$ classical spectral curve \eqref{SpecCurveGenus} can be expressed completely in terms of the ones produced by the topological recursion on the classical spectral curve of genus $0$ given by \eqref{SpecCurveGenus0}. Note that the relation between both sets is non-trivial and corresponds to matching \eqref{GeneralForm2R1} with \eqref{MainTheoAsymp} at each order in $N^{-k}$ in the expansion. To the knowledge of the author, no known other examples of such correspondence are known. Thus, it would be interesting to see if one can obtain such results directly from the topological recursion formalism. In particular, one could expect some kinds of factorization of the recursion kernel in such a way that this correspondence can be tracked down at each step of the recursion. This could provide a general technique to take benefit from symmetries in the topological recursion formalism as initiated in \cite{DuninNorburyOrantinSpectralCurveWithSymmetries2019} in a different context. 

Finally, the results are also interesting in the context of quantization of classical spectral curves using the topological recursion. Indeed, the quantization method proposed in \cite{Quant1,Quant2,Quant3} in relation with Lax pairs and integrable systems implies the use of formal $\frac{1}{N}$-transseries when the genus of the classical spectral curve in non-zero (when the genus is vanishing, it reduces to standard formal power series or WKB series). These $\frac{1}{N}$-transseries include oscillatory terms and the construction is still at the formal level. But current works \cite{Costin2,Costin3,Nikolaev1,Nikolaev2,Nikolaev3} suggests that one may resum these $\frac{1}{N}$-transseries to obtain analytic objects and thus recover rigorously equalities similar to \eqref{MainTheoAsymp} in a general way. In this context, the present work could a non-trivial check of arbitrary genus for these resummation techniques and shed some light on how to deal properly with $\frac{1}{N}$-transseries.  

\subsection{Possible generalizations}

It is natural to wonder if the results presented in this article may be generalized or extended in other ways. We list here some of the possibilities.
\begin{itemize}\item \underline{$m$ depending on $N$}: The large $N$ asymptotic expansion given in Theorem \ref{MainTheo} is performed for a fixed value of $m$. It is natural to wonder if such asymptotic expansions could be extended to situations when $m$ depends on $N$, i.e. sequences of integers $(m_N)_{N\geq 1}$. If the sequence of integers converges, then there exists $N_0\in \mathbb{N}$ such that for all $N\geq N_0$, $m_N=M\in \mathbb{N}$. Hence, the results presented in this article trivially apply for $N\geq N_0$ and so does the large $N$ asymptotic expansion of Theorem \ref{MainTheo} with $m$ replaced by $M$. Other regimes when $m_N\to +\infty$ or $m_N$ periodic could be studied. However, we recall here that when $m_N>N$ then the Toeplitz determinant trivially equals $\epsilon^N$. Interesting sequences are for example strictly increasing sequences $(m_N)_{N\geq 1}$ with $m_N=o(N)$. In such situations, the factorization of Proposition \ref{Factorization} still holds with $N=n_1(N) m_N+n_2(N)$  with $n_2(N)<m_N$ at each step (we added the dependence in $N$ of $n_1$ and $n_2$ for clarity). Thus, we get
\beqq   \ln D_N(m,\epsilon)=(m_N-n_2(N))\ln D_{n_1(N)}(\mathds{1}_{[-\pi \epsilon,\pi \epsilon]})+n_2(N)\ln D_{n_1(N)+1}(\mathds{1}_{[-\pi \epsilon,\pi \epsilon]})\eeqq
Note that if $m_N=o(N)$ then $n_1(N)\overset{n\to \infty}{\to} +\infty$ so that the asymptotic expansion of both $\ln D_{n_1(N)}(\mathds{1}_{[-\pi \epsilon,\pi \epsilon]})$ and $\ln D_{n_1(N)+1}(\mathds{1}_{[-\pi \epsilon,\pi \epsilon]})$ is still given by Proposition \ref{OneCutCase} with $n_1(N)$ (resp. $n_1(N)+1$) replacing $n$. One finally need to multiply these expansions by $m_N-n_2(N)$ (resp. $n_2(N)$) and order the asymptotic expansion properly. However, the main difficulty is that $(n_2(N))_{N\geq 1}$ is no longer necessary bounded. Indeed, the sequence $(n_2(N))_{N\geq 1}$ may or may not be bounded or periodic depending on the choice of $(m_N)_{N\geq 1}$. Regarding the situation, one has to choose some appropriate sets of parameters to express the asymptotic expansion nicely, but the present strategy presented in this article still holds. 
\item \underline{Non-integer powers of unitary matrices}: Another natural generalization is to consider real powers of unitary matrices: $U_N^{\, t}$ with $t\in \mathbb{R}_+$ and consider the probability that all its eigenvalues are located in $\mathcal{C}_1(\epsilon)$as in \cite{UnitaryMarchal2014}. This question is equivalent to ask for the eigenvalues of $U_N$ to be located in a union of arc-intervals (the exact expression is given in equation $4.3$ of \cite{UnitaryMarchal2014}). However, unlike the present case, the discrete rotational symmetry does not apply so the factorization of the Toeplitz determinant fails. Nevertheless, the asymptotic expansion of Theorem \ref{MainTheo} may be naively generalized by setting  $m=t$, $n_1=\left\lfloor \frac{N}{t}\right\rfloor$ and $n_2=N-\left\lfloor \frac{N}{t}\right\rfloor t$. It would be interesting to study if the first orders of the asymptotic expansion obtained with this adaptation are still correct.   
\item  \underline{$\beta$-CUE generalization}: Finally one may consider the same problem but with matrices drawn from the $\beta$-CUE, i.e. $D_N(m,\epsilon,\beta)=\frac{1}{(2\pi)^N\,N!}\int_{I_m(\epsilon)^{\, N}} d\theta_1\dots d\theta_N \underset{1\leq p<q\leq N}{\prod}\left|e^{i \theta_p}-e^{i\theta_q}\right|^{2\beta}$. In this case, an asymptotic expansion similar to Proposition \ref{TopRecGenus} is known \cite{BorotGuionnetKozBetaMultiCut2015}. However, this case cannot be written in terms of Toeplitz determinants and thus the present strategy cannot be carried out. The only possible cases for which one may hope some simple generalizations of the present strategy are $\beta\in\left\{\frac{1}{2},2\right\}$ because connections with bi-orthogonal polynomials and Pfaffian exist. However, on the topological recursion side, $\beta\neq 1$ implies drastic changes \cite{MarchalBetaEnsembleZhukovsky2011,MarchalChekovEynardBetaTopRecSectorwise2011} in the recursion even for genus $0$ spectral curves and much less is known in this situation. 
\end{itemize}

\section*{Acknowledgments}
I would like to thank C. Charlier for suggesting the reference to B. Fahs PhD. thesis and several other references. I would like to thanks B. Fahs for fruitful discussions and N. Orantin for suggesting reference \cite{DuninNorburyOrantinSpectralCurveWithSymmetries2019}.

\section*{Statements and declarations}
The author declares that there is no conflict of interest.

\bibliographystyle{plain}
\bibliography{ToeplitzMultiIntervalWithSymmetry}

\begin{appendix}

\section{Topological recursion details in the one-cut case \label{AppendixComputationTopRec}}
In this section, we list the intermediate steps obtained in the computation of the topological recursion of \cite{EynardOrantin2007} applied to the classical spectral curve of genus $0$ given by \eqref{SpecCurveGenus0} that may be parametrized globally as:
\bea x(z)&=&\frac{1}{2}\tan \left(\frac{\pi \epsilon}{2}\right) \left(z+\frac{1}{z}\right)\cr
y(z)&=&\frac{8 z^3}{\sin \left(\frac{\pi \epsilon}{2}\right) (z-1)(z+1) \left( \left(\tan^2 \frac{\pi \epsilon}{2}\right) z^4+2\left( \tan^2 \frac{\pi \epsilon}{2}+2\right)z^2+  \tan^2 \frac{\pi \epsilon}{2}\right)} 
\eea
with the global involution $\bar{z}=\frac{1}{z}$. The parametrization implies that there are two simple branchpoints at $z=1$ and $z=-1$. The definition of the topological recursion in this particular case is presented in Appendix B of \cite{ToeplitzMarchal2020}. It provides the following intermediate steps:
\bea \omega_2^{(0)}(z_1,z_2)&=&\frac{dz_1 \, dz_2}{(z_1-z_2)^2}\cr
\omega_{1}^{(1)}(z)&=& \frac{z }{2(z-1)^2(z+1)^2 \cos \frac{\pi \epsilon}{2} }\cr
\omega_{3}^{(0)}(z)&=&0\cr
\omega_{4}^{(0)}(z)&=&0\cr
\omega_{2}^{(1)}(z_1,z_2)&=&\frac{z_1^2z_2^2+z_1^2+z_2^2+4z_1z_2+1}{4 \left(\cos\frac{\pi \epsilon}{2}\right)^2 (z_1-1)^2(z_1+1)^2(z_2-1)^2(z_2+1)^2} \cr
\omega_1^{(2)}(z)&=&\frac{z\left(\left(13\tan^2 \frac{\pi \epsilon}{2} +4\right)z^4+\left(10\tan^2 \frac{\pi \epsilon}{2} +28\right)z^2+13\tan^2 \frac{\pi \epsilon}{2} +4\right)}{32 \left(\cos\frac{\pi \epsilon}{2}\right) (z+1)^4(z-1)^4}\cr
\omega_3^{(1)}(z_1,z_2,z_3)&=&\frac{(z_1z_2+z_1z_3+z_2z_3+1)(z_1z_2z_3+z_1+z_2+z_3)}{\left(\cos\frac{\pi \epsilon}{2}\right)^3(z_1^2-1)^2(z_2^2-1)^2(z_3^2-1)^2}\cr
\omega_1^{(3)}(z)&=&\frac{z}{\left(\cos\frac{\pi \epsilon}{2}\right)(z^2-1)^6}\Big(\left(413 \tan^4 \frac{\pi \epsilon}{2}+268 \tan^2 \frac{\pi \epsilon}{2}+8 \right)(z^8+1)\cr
&& +\left(580 \tan^4 \frac{\pi \epsilon}{2}+1880 \tan^2 \frac{\pi \epsilon}{2}+688\right)(z^6+z^2)\cr
&&  +\left(1614 \tan^4 \frac{\pi \epsilon}{2}+2904 \tan^2 \frac{\pi \epsilon}{2}+2208\right)z^4\Big)\cr
&&
\eea

\end{appendix}

\end{document}